\documentclass[final,12pt,3p]{elsarticle}

\usepackage[T1]{fontenc}
\usepackage[utf8]{inputenc}
\usepackage{amsmath,amssymb,amsthm,mathtools}
\usepackage{enumitem}
\usepackage{geometry}
\usepackage{verbatim} 

\usepackage{lineno}
\usepackage[expansion=false]{microtype}
\usepackage{xcolor}
\usepackage{hyperref}

\usepackage{caption}
\usepackage{adjustbox}
\usepackage{diagbox}

\usepackage{algorithm}
\usepackage{algorithmic}
\usepackage{tikz}
\usetikzlibrary{shapes}
\usepackage{subcaption}
\usepackage{mathrsfs} 
\usepackage{multirow}

\biboptions{sort&compress}
\hypersetup{
  colorlinks=true,
  linkcolor=blue,
  citecolor=blue,
  urlcolor=blue,
  pdftitle={Convex Structure-Preserving EDMD: Passive Koopman Matrices as Q-Contractions},
  pdfauthor={Matthias Ehrhardt and Achraf Zinihi},
}
\allowdisplaybreaks 

\theoremstyle{plain}
\newtheorem{theorem}{Theorem}[section]
\newtheorem{proposition}[theorem]{Proposition}

\newtheorem{corollary}[theorem]{Corollary}

\theoremstyle{definition}
\newtheorem{definition}[theorem]{Definition}

\theoremstyle{remark}

\newtheorem{example}[theorem]{Example}

\newcommand{\R}{\mathbb{R}}
\newcommand{\E}{\mathbb{E}}

\newcommand{\Hamil}{H}
\newcommand{\grad}{\nabla}

\newcommand{\Lift}{\Psi}

\newcommand{\Xspace}{\mathcal{X}}

\newcommand{\SpH}{\mathcal{S}_{\mathrm{pH}}}
\newcommand{\JK}{J_{\mathrm{K}}}
\newcommand{\RK}{R_{\mathrm{K}}}
\newcommand{\HK}{\Hamil_{\mathrm{K}}}
\newcommand{\yK}{y_{\mathrm{K}}}

\DeclareMathOperator{\diag}{diag}
\DeclareMathOperator{\dist}{dist}

\DeclareMathOperator{\tr}{tr} 
\DeclareMathOperator{\ran}{range}
\DeclareMathOperator{\Cay}{Cay} 

\usepackage{graphicx}
\usepackage{booktabs}

\usepackage[figcolor=white]{todonotes}

\begin{document}

\begin{frontmatter}

\title{Convex Structure-Preserving EDMD:\\ Passive Koopman Matrices as $Q$-Contractions}

\author[BUW]{Matthias Ehrhardt\corref{Corr}}
\cortext[Corr]{Corresponding author}
\ead{ehrhardt@uni-wuppertal.de}

\author[BUW]{Achraf Zinihi}
\ead{achrafzinihi99@gmail.com}

\address[BUW]{Applied and Computational Mathematics, University of Wuppertal,\\
Gaußstrasse 20, Wuppertal, 42119, Germany}


\begin{abstract}
Extended dynamic mode decomposition (EDMD) approximates the Koopman operator using a matrix based on a finite set of observables. However, the least-squares matrix of a port-Hamiltonian (pH) system is generally neither passive nor invariant-preserving.
We demonstrate that the implicit midpoint discretization of a lifted pH model precisely preserves energy balance, passivity, and all linear Casimirs for every step size.
For a fixed storage metric $Q$, the resulting Koopman matrices are the Cayley transforms of $(J-R)Q$, 
where $J$ is skew-symmetric and $R$ is positive semidefinite. This set does not depend on the step size and is the closure of all Q-contractions.
This set is LMI-representable; therefore, structure-preserving EDMD, including ports through a discrete-time Kalman-Yakubovich-Popov condition and Casimirs 
through affine constraints, is a semidefinite program.
The constrained estimate is a projection of the unconstrained estimate in the data-weighted norm and cannot increase estimation error when the population target is pH-compatible.
For every $Q$-contraction, a certified bound on the lifted prediction error grows linearly at most in time, 
whereas it grows geometrically for unconstrained models with Q-gain above one.
Benchmarks on a damped pendulum and a damped rigid body demonstrate that constrained models never generate energy and that imposed Casimirs are maintained with solver accuracy (approximately $10^{-8}$). 
In contrast, EDMD may diverge or lose the Casimir with noisy data. 
The benchmarks also demonstrate that when $Q$ is the physical energy, the constrained models are less accurate than EDMD and can exhibit significant structural bias.
\end{abstract}

\begin{keyword}
Koopman operator \sep extended dynamic mode decomposition \sep port-Hamiltonian systems \sep passivity \sep Cayley transform \sep linear matrix inequalities \sep Casimir functions

\vspace{.2cm}
\end{keyword}

\journal{Systems \& Control Letters}

\end{frontmatter}


\section{Introduction}\label{S1} 
Extended dynamic mode decomposition (EDMD) \cite{Williams2015, Korda2018JNLS} approximates the Koopman operator \cite{Koopman1931, Mezic2005, Brunton2022, Zinihi2026KoopmanPINN} of a nonlinear system with a matrix that acts on a finite dictionary of observables.
It has become a standard tool for prediction and predictive control \cite{Korda2018Automatica, Bruder2021, Zinihi2026KoopmanSIRSD}.
However, for port-Hamiltonian (pH) systems \cite{vanderSchaft2014, vanderSchaft2017} the unconstrained least-squares matrix is structurally agnostic: 
it need not be passive with respect to any storage function, nor does it need to conserve the Casimirs of the system. Long-horizon simulations often exhibit spurious energy drift.

Structure-preserving fits are available for linear pH systems in state coordinates
\cite{Morandin2023, Cherifi2022, Zwerschke2026}.,
The operator-theoretic structure of the Koopman generator of nonlinear pH systems has also been clarified \cite{Preciado2026}.
In this letter, we address the discrete-time identification problem in lifted coordinates and demonstrate that it is a convex program for a fixed storage metric.


The main contributions of this letter are as follows:
First, we demonstrate that the implicit midpoint discretization of a lifted pH model exactly preserves the energy balance, passivity, and linear Casimirs for all step sizes. Second, we show that the resulting Koopman matrices with the storage function $\tfrac{1}{2}z^\top Qz$ are precisely the $Q$-contractions, 
independent of the step size $\Delta t$.
Consequently, structure-preserving EDMD with ports and Casimirs is a semidefinite program (SDP).
The solution to this program is a projection of the unconstrained estimate in the data-weighted norm. 
Therefore, the constraint cannot increase the estimation error when the population target is pH-compatible. 
Furthermore, the lifted prediction error of the model grows at most linearly with time.
Three benchmarks quantify the benefits and costs of the constraints.


The Cayley (bilinear, Tustin) transform and the discrete-time Kalman-Yakubovich-Popov (KYP) Lemma are classical. 
Discrete-time pH systems based on symplectic integration and collocation have been studied in references \cite{Talasila2006, Kotyczka2018, Kotyczka2019, Mehrmann2019, CelledoniHoiseth2017}.
This letter contributes their combination in the lifted Koopman setting, yielding an identification method that is exactly passive and Casimir-preserving while remaining convex.
A systematic choice of the storage metric, design of admissible dictionaries, and sharper error bounds are topics for future research.


Three recent works are the most similar to ours. 
Lortie et al.~\cite{Lortie2024} and Mamakoukas et al.~\cite{Mamakoukas2020} constrain the EDMD matrix to be stable;
however, stability alone does not prevent transient energy growth. 
Singh et al.~\cite{Singh2026}  learn a lifting with a neural network together with a pH generator, discretize it using a Cayley (midpoint) step, and prove zero-input contractivity.
Thus, the midpoint step itself is not novel. 
This letter adds an exact description of the entire passive set independent of $\Delta t$, a direct convex fit of $A_d$ from snapshot pairs with ports and Casimirs, a projection property, and a quantitative account of accuracy costs.


\section{The Problem Setting}\label{S2}
We consider input-state-output pH systems of the form
\begin{equation}\label{E1}
   \dot{x} = \bigl[ J(x)-R(x) \bigr] \grad\Hamil(x)+G(x)u, \quad\text{with}\quad
   y = G(x)^\top\grad\Hamil(x),
\end{equation}
on $\Xspace\subseteq\R^n$, with a \textit{Hamiltonian} $\Hamil\in C^2(\Xspace)$ that is bounded below, with a skew-symmetric interconnection matrix $J(x)=-J(x)^\top$, and a symmetric positive semidefinite dissipation matrix $R(x)=R(x)^\top \succeq 0$.
Along solutions of \eqref{E1}, we have the \textit{dissipation inequality}
\begin{equation}\label{eq:dissipation}
  \frac{d}{d t}\Hamil(x(t)) = -\grad\Hamil^\top R\grad\Hamil+u^\top y\le u^\top y,
\end{equation}
i.e., the pH system \eqref{E1} is passive with storage function $\Hamil$ \cite{Willems1972}.
We recall that a function $C$ with the property $J\grad C\equiv 0$ is a Casimir \cite{MarsdenRatiu1999}.

Next, let $\Lift = (\psi_1,\dots,\psi_N)^\top$ be a \textit{dictionary of observables} and $z=\Lift(x)$ the \textit{lifted state}.
Given snapshot pairs $(x_k,x_k^+)$, $k=1,\dots,K$, where $x_k^+$ denotes the state one sampling period $\Delta t$ after $x_k$, standard EDMD computes
\begin{equation}\label{E2}
  \widehat{K} = \arg\min_{K\in\R^{N\times N}} \bigl\lVert Z_Y - K Z_X \bigr\rVert_F^2,
\end{equation}
where
\begin{equation*}
   Z_X = \bigl[ \Lift(x_1),\dots,\Lift(x_K) \bigr], \quad\text{and}\quad 
   Z_Y = \bigl[ \Lift(x_1^+),\dots,\Lift(x_K^+) \bigr].
\end{equation*}
For the continuous-time model class, we use the \textit{lifted pH system}
\begin{equation}\label{E3}
    \dot{z} = \Bigl( \JK - \RK + \sum_{j=1}^m u_j S_j \Bigr)Qz + Bu, \quad\text{with} \quad
    y_{\mathrm{K},j} = (Qz)^\top S_j(Qz) +\bigl(B^\top Qz\bigr)_j,
\end{equation}
with $Q = Q^\top\succ 0$, $\JK = -\JK^\top$, $\RK = \RK^\top\succeq 0$, arbitrary $S_j\in\R^{N\times N}$, $B\in\R^{N\times m}$, and lifted storage function $\HK(z) = \tfrac{1}{2}z^\top Qz$.
%
Direct computation using $Qz^TJK(Qz)=0$ yields the following 
(cf.\ \eqref{eq:dissipation})
\begin{equation*}
     \frac{d}{d t}\HK = -(Qz)^\top \RK(Qz) + \yK^\top u \le \yK^\top u,
\end{equation*} 
for all parameter values. 
Thus, the lifted pH system~\eqref{E3} is passive by construction.
The bilinear terms $S_j$ are necessary for an exact representation of the energy balance of control-affine systems.
Throughout this paper, the storage metric $Q$ is given.
Section~\ref{S6} discusses how it can be chosen.

\section{Exact Discrete-Time Structure Preservation}\label{S3}
Continuous-time passivity \eqref{eq:dissipation} is not inherited by arbitrary discretizations.
For example, the energy increment of the explicit Euler method contains the nonnegative remainder $\tfrac{\Delta t^2}{2}f^\top Qf$, where $f$ is the right-hand side of \eqref{E3}, which can destroy passivity.
Therefore, we work with the implicit midpoint rule.

\begin{definition}[Discrete-time lifted pH system]\label{D3.1}
Given the parameters of \eqref{E3} and $\Delta t>0$, the 
\textit{implicit midpoint scheme} reads
\begin{equation}\label{E4}
  z_{k+1} = z_k + \Delta t \Bigl[\Bigl(\JK-\RK + \sum_{j=1}^m  u_{k,j} S_j \Bigr) Q z_{k+1/2} + B u_k \Bigr], \quad\text{where}\quad 
  z_{k+1/2} := \tfrac{1}{2}(z_k+z_{k+1}),
\end{equation}
and the outputs are $y_{\mathrm{K},k,j}:=(Qz_{k+1/2})^\top S_j(Qz_{k+1/2})+(B^\top Qz_{k+1/2})_j$.
\end{definition}

The scheme requires solving one linear system per step, with the following matrix $I-\tfrac{\Delta t}{2}(\JK-\RK+\sum_ju_{k,j}S_j)Q$.
For $S_j=0$, this matrix is invertible for all $\Delta t>0$ (Proposition~\ref{P4.1}). 
For $S_j\neq 0$, unique solvability depends on the values of both the time step and the input amplitude.

\begin{theorem}[Exact discrete structure preservation]\label{T3.2}
Let $\Delta t>0$, and suppose that \eqref{E4} is uniquely solvable at every step.
Then the scheme \eqref{E4} has the following properties:
\begin{enumerate}
\item[$(i)$] \emph{exact energy balance} (as a discrete analogue to the dissipation inequality \eqref{eq:dissipation}):
\begin{equation}\label{E5}
    \HK(z_{k+1}) - \HK(z_k) = -\Delta t\,(Qz_{k+1/2})^\top\RK(Qz_{k+1/2})+\Delta t\,u_k^\top y_{\mathrm{K},k} \le\Delta t\,u_k^\top y_{\mathrm{K},k}.
\end{equation}
Hence, discrete-time passivity with storage function $\HK$;
\item[$(ii)$] \emph{exact Casimir preservation:} 
if $c\in\ker\JK\cap\ker\RK$ and, in the forced case, $c^\top S_j=0$ and $c^\top B=0$, then $c^\top z_{k+1}=c^\top z_k$;
\item[$(iii)$] \emph{second-order accuracy:} 
the scheme is a second-order approximation of the flow of \eqref{E3}.
\end{enumerate}
\end{theorem}
\begin{proof}
$(i)$ Since $Q=Q^\top$, we have the energy balance $\HK(z_{k+1})-\HK(z_k) = (z_{k+1}-z_k)^\top Qz_{k+1/2}$.
Substituting the implicit midpoint scheme \eqref{E4} and using $(Qz_{k+1/2})^\top \JK(Qz_{k+1/2}) = 0$ gives \eqref{E5}.

$(ii)$ Multiplying the implicit midpoint scheme \eqref{E4} by $c^\top$ and using $c^\top \JK = c^\top \RK = 0$, which follows from $\JK c=\RK c=0$ by skew symmetry and symmetry, respectively, together with $c^\top S_j=0$ and $c^\top B=0$, gives $c^\top(z_{k+1}-z_k)=0$.

$(iii)$ This is the classical order of accuracy of the implicit midpoint rule.
\end{proof}

Theorem~\ref{T3.2} shows that the continuous pH structure carries over to a discrete energy inequality \textit{exactly}, 
at the cost of one linear solve per step, and that the Casimirs are carried along.
For a nonlinear lifted Hamiltonian, 
one replaces $Qz_{k+1/2}$ with a discrete gradient $\overline{\grad} \HK(z_k,z_{k+1})$, cf.~\cite{Gonzalez1996, Kotyczka2019, CelledoniHoiseth2017}, and the above mentioned properties $(i)$--$(ii)$ remain valid verbatim.

\section{Passive Koopman Matrices are \texorpdfstring{$Q$}{Q}-Contractions}\label{S4}
The structure-preserving version of \eqref{E2} restricts $K$ to the set of Koopman matrices that admit a pH realization. The following result identifies this set.

\begin{proposition}[Cayley correspondence]\label{P4.1}
Let $Q\succ 0$ and $\Delta t>0$.
For $\mathcal{A}\in\R^{N\times N}$, the following are equivalent:
\begin{enumerate}
\item[$(a)$] $\mathcal{A}=(\JK-\RK)Q$ for some $\JK=-\JK^\top$ and $\RK\succeq 0$;
\item[$(b)$] $Q\mathcal{A}+\mathcal{A}^\top Q\preceq 0$.
\end{enumerate}
In this case, $I-\tfrac{\Delta t}{2}\mathcal{A}$ is invertible, and the \emph{Cayley transform}
\begin{equation*}
     A_d = \Cay_{\Delta t}(\mathcal{A}) := \Bigl(I-\frac{\Delta t}{2}\mathcal{A}\Bigr)^{-1}\Bigl(I+\frac{\Delta t}{2}\mathcal{A}\Bigr)
\end{equation*}
satisfies the \emph{$Q$-contraction property} $A_d^\top QA_d\preceq Q$ and inherits the linear invariants: 
$c^\top\mathcal{A}=0$ implies $c^\top A_d=c^\top$.
Conversely, every $A_d$ with $A_d^\top QA_d\preceq Q$ and $-1\notin\sigma(A_d)$ is the Cayley transform of a unique such $\mathcal{A}$.
\end{proposition}
\begin{proof}
$(a) \Rightarrow (b)$: 
$Q\mathcal{A}+\mathcal{A}^\top Q = -2Q\RK Q\preceq 0$.

$(b) \Rightarrow (a)$: 
Take $\JK = \tfrac{1}{2}(\mathcal{A} Q^{-1} - Q^{-1} \mathcal{A}^\top)$ and $\RK = -\tfrac{1}{2}(\mathcal{A}Q^{-1} + Q^{-1} \mathcal{A}^\top)$. 
$\RK\succeq 0$ is exactly $(b)$ after congruence with $Q^{-1}$.
By $(b)$, the spectrum of $\mathcal{A}$ lies in the closed left half-plane, so $2/\Delta t$ is not an eigenvalue and $I-\tfrac{\Delta t}{2} \mathcal{A}$ is invertible.
With $W = I-\tfrac{\Delta t}{2} \mathcal{A}$ and $V=I+\tfrac{\Delta t}{2} \mathcal{A}$, which commute,
\begin{equation*}
  A_d^\top Q A_d - Q = W^{-\top} \bigl( V^\top Q V - W^\top Q W \bigr) W^{-1},
\quad \text{and}\quad
V^\top QV-W^\top QW=\Delta t\,(\mathcal{A}^\top Q+Q\mathcal{A})\preceq 0.
\end{equation*}
If $\mathcal{A}^\top c=0$, then $W^\top c=c$, hence $W^{-\top}c=c$ and $A_d^\top c=V^\top W^{-\top}c=V^\top c=c$.
For the converse, the inverse transform is $\mathcal{A} = \tfrac{2}{\Delta t}(A_d+I)^{-1}(A_d-I)$, and the computation above, read backwards, shows that $A_d^\top QA_d\preceq Q$ implies $(b)$.
\end{proof}

We write $\lVert v\rVert_Q := (v^\top Q v)^{1/2}$, and use the same symbol for the induced matrix norm, so that $A_d^\top Q A_d \preceq Q$ is equivalent to $\lVert A_d \rVert_Q\le 1$. We denote the set of $Q$-contractions by $\SpH(Q) := \{A_d\in\R^{N\times N}\colon A_d^\top Q A_d\preceq Q\}$.
The subset $\{A_d \in \SpH(Q)\colon -1\notin \sigma(A_d)\}$ in the converse part of Proposition~\ref{P4.1} does not involve $\Delta t$; 
only the generator $\mathcal{A}$ attached to a given $A_d$ does. 
Furthermore, if $A_d\in\SpH(Q)$, then $A_\varepsilon=(1-\varepsilon)A_d$ satisfies $\lVert A_\varepsilon\rVert_Q\le 1-\varepsilon$, so $-1\notin\sigma(A_\varepsilon)$ and $A_\varepsilon\to A_d$ as $\varepsilon\to 0$. 
Therefore, the Cayley images of $Q$-dissipative generators are dense in the closed set $\SpH(Q)$. 
Therefore, searching over all $Q$-contractions loses nothing, up to closure, compared with searching over Cayley-discretized pH models, and no derivative data or generator parametrization is needed.

\begin{corollary}[Convex structure-preserving EDMD]\label{C4.3}
For a fixed $Q\succ 0$, the set $\SpH(Q)$ is closed, convex, and representable by a linear matrix inequality (LMI)
\begin{equation}\label{E6}
\begin{bmatrix}
Q & A_d^\top Q\\  QA_d & Q
\end{bmatrix}
\succeq 0,
\end{equation}
so that $\min\limits_{A_d\in\SpH(Q)} \lVert Z_Y-A_dZ_X\rVert_F^2$ is a convex SDP. 
With ports, the passivity of $(A_d, B_d, C_d, D_d)$ with storage function $\tfrac{1}{2}z^\top Qz$ and supply rate $u^\top y$ is equivalent to the discrete-time KYP condition
\begin{equation}\label{E7}
  \begin{bmatrix}
    Q & C_d^\top & A_d^\top\\
    C_d & D_d+D_d^\top & B_d^\top\\
    A_d & B_d & Q^{-1}
  \end{bmatrix}\succeq 0,
\end{equation}
which is again an LMI in $(A_d, B_d, C_d, D_d)$ for a fixed $Q$.
Imposing the linear Casimir $C_0^\top z$ adds the affine constraints $C_0^\top A_d = C_0^\top$ and $C_0^\top B_d=0$, which preserve convexity.
\end{corollary}
\begin{proof}
Both LMIs follow from a Schur complement applied to
\begin{equation*}
   [A_d\ B_d]^\top Q[A_d\ B_d] - 
\begin{bmatrix}
  Q&C_d^\top\\ C_d&D_d+D_d^\top
\end{bmatrix}
\preceq 0,
\end{equation*}
which is the dissipation inequality $\tfrac{1}{2}z_{k+1}^\top Q z_{k+1} - \tfrac{1}{2}z_k^\top Q z_k \le u_k^\top y_k$ written as a quadratic form in $(z_k, u_k)$. 
For \eqref{E6}, set $B_d=0$, $C_d=0$, and $D_d=0$, delete the resulting zero block row and column, and multiply the last block row and column by $Q$.
\end{proof}

The matrices $(A_d,B_d)$ above describe models that are linear in the input. 
The bilinear terms $S_j$ in \eqref{E3} cause the midpoint map to be nonlinear in $(z_k,u_k)$, so Corollary~\ref{C4.3} covers the case $S_j=0$, which is also the case used in Section~\ref{S7}. 
The Casimir constraints never render the program infeasible.  
For $C_0$ with full column rank, $A_d = Q^{-1} C_0 (C_0^\top Q^{-1} C_0)^{-1} C_0^\top$, $B_d = 0$ is the $Q$-orthogonal projector onto $\ran(Q^{-1}C_0)$. 
Thus, it is a $Q$-contraction, that satisfies both constraints.

Corollary~\ref{C4.3} is the algorithmic core of the approach. 
For a fixed storage metric, structure-preserving EDMD is a convex program, requiring only snapshot pairs, not derivative data.
Its continuous-time counterpart is the LMI $Q\mathcal{A}+\mathcal{A}^\top Q\preceq 0$ of Proposition~\ref{P4.1}$(b)$.
The nonconvexity of the identification problem lies entirely in jointly estimating $Q$ and the remaining parameters.
This is why $Q$ is fixed before the regression.
After solving the SDP, a pH realization $(\JK, \RK)$ with $\Cay_{\Delta t} ( (\JK-\RK) Q) =  A_d$ is obtained using  the inverse Cayley transform and the formulas in the proof of Proposition~\ref{P4.1}, provided that $-1\notin\sigma(A_d)$.


\section{Projection Property and Error Growth}\label{S5}
Given $z_k = \Lift(x_k)$ and $z_k^+ = \Lift(x_k^+)$, we can write the regression equation \eqref{E2}, possibly augmented by an input matrix, as a linear least-squares problem in a parameter vector $\xi\in\R^p$: 
for $\xi = \operatorname{vec}(A_d)$, we have $A_dz_k = \Phi(z_k)\xi$ with $\Phi(z_k) = z_k^\top \otimes I_N$.

Let $\mathcal{C}\subseteq\R^p$ be the closed convex constraint set of Corollary~\ref{C4.3}, and define the \emph{population target}, cf.\ \cite[Section~5.3.2]{chen2026beyond}, (also called $L_2(\nu)$-projection in the foundational paper \cite{Korda2018JNLS}):
%
%
%
\begin{equation*}
   \xi^{\star} = \arg\min_{\xi\in\R^p}\ \E \bigl\lVert z^+ -\Phi(z)\xi\bigr\rVert^2,
\end{equation*}
where the expectation is taken with respect to the sampling distribution of the snapshot pairs.

Define the \textit{empirical Gram matrix} $\widehat{G} := \frac{1}{K}\sum_{k=1}^K \Phi(z_k)^\top \Phi(z_k) = \bigl( \tfrac{1}{K} Z_X Z_X^\top \bigr)\otimes I_N$, which is assumed to be positive definite.
Define the $\widehat{G}$-induced norm as follows: $\lVert \xi \rVert_{\widehat{G}} := ( \xi^\top \widehat{G} \xi )^{1/2}$.

\begin{proposition}[Nonexpansiveness of the structural constraint]\label{P5.1}
Let $\mathcal{C}$ be a closed, convex set. 
Let $\widehat{\xi}$ be the unconstrained least-squares estimate, 
$\widehat{\xi}_{\mathcal{C}}$ the constrained one. 
Let the population target be denoted by $\xi^{\star}$ 
Then $\widehat{\xi}_{\mathcal{C}} = \Pi_{\mathcal{C}}(\widehat{\xi})$ is the metric projection of $\widehat{\xi}$ onto $\mathcal{C}$ in the norm $\lVert\cdot\rVert_{\widehat{G}}$, and
\begin{equation*}
  \lVert \widehat{\xi}_{\mathcal{C}}-\xi^{\star} \rVert_{\widehat{G}}
  \le \lVert \widehat{\xi} - \xi^{\star} \rVert_{\widehat{G}}
  + \dist_{\widehat{G}}(\xi^{\star},\mathcal{C}).
\end{equation*}
In particular, if $\xi^{\star}\in\mathcal{C}$, that is, if the population target is itself pH-compatible, then the constrained estimator is never worse than the unconstrained one:
$\lVert \widehat{\xi}_{\mathcal{C}} -   
   \xi^{\star} \rVert_{\widehat{G}} 
  \le \lVert \widehat{\xi} - \xi^{\star} \rVert_{\widehat{G}}$.
\end{proposition}
\begin{proof}
By the normal equations, the empirical loss satisfies
\begin{equation*}
  \frac{1}{K} \sum_{k=1}^K \bigl\lVert z_k^+ - \Phi(z_k) \xi \bigr\rVert^2
 = \frac{1}{K} \sum_{k=1}^K \bigl\lVert z_k^+ -\Phi(z_k) \widehat{\xi} \bigr\rVert^2
  +\lVert \xi - \widehat{\xi} \rVert_{\widehat{G}}^2,
\end{equation*}
so minimizing it over $\mathcal{C}$ is the same as projecting $\widehat{\xi}$ onto $\mathcal{C}$ in the metric $\lVert \cdot \rVert_{\widehat{G}}$.
The metric projection onto a closed convex set in a Hilbert space is $1$-Lipschitz.
Now, using $\bar{\xi} := \Pi_{\mathcal{C}}(\xi^{\star})$, we obtain
\begin{equation*}
  \lVert \Pi_{\mathcal{C}}(\widehat{\xi})-\xi^{\star} \rVert_{\widehat{G}}
  \le \lVert \Pi_{\mathcal{C}}(\widehat{\xi})-\bar{\xi} \rVert_{\widehat{G}}
  + \lVert \bar{\xi} - \xi^{\star} \rVert_{\widehat{G}}
  \le \lVert \widehat{\xi} - \xi^{\star} \rVert_{\widehat{G}} + \dist_{\widehat{G}}(\xi^{\star},\mathcal{C}).
\end{equation*}
If $\xi^{\star}\in\mathcal{C}$, then $\bar{\xi}=\xi^{\star}$ and the distance term vanishes.
\end{proof}

Because $\widehat{\xi}_{\mathcal{C}}$ is a projection, the obtuse-angle property $\langle \widehat{\xi} - \widehat{\xi}_{\mathcal{C}},\,\zeta - \widehat{\xi}_{\mathcal{C}} \rangle_{\widehat{G}} \le 0$ holds for every $\zeta \in \mathcal{C}$, which gives the Pythagorean inequality
\begin{equation*}
  \lVert \widehat{\xi}_{\mathcal{C}}-\zeta \rVert_{\widehat{G}}^2 + \lVert \widehat{\xi}-\widehat{\xi}_{\mathcal{C}} \rVert_{\widehat{G}}^2
  \le \lVert \widehat{\xi} - \zeta \rVert_{\widehat{G}}^2 \qquad \forall\,\zeta\in\mathcal{C} .
\end{equation*}
Thus the constrained estimate is never farther than EDMD from any passive reference model, and it is strictly closer unless $\widehat{\xi}\in\mathcal{C}$. This is a deterministic statement and needs no assumption on $\xi^\star$.

The norm in Proposition~\ref{P5.1} is the data-weighted norm induced by the regression, not the Frobenius norm of the parameter matrices. 
This statement applies to the pure least-squares problem. Additional penalties, such as those imposed for  output mismatch, modify the projection.
The term $\dist_{\widehat{G}}(\xi^{\star},\mathcal{C})$ is the \textit{structural bias}.
Proposition~\ref{P5.1} thus reduces the statistical question to a deterministic approximation-theoretic one, namely the size of the structural bias. 
This proposition is consistent with the robustness of constrained fits in the small-sample regime, because 
the variance contribution is not amplified, whereas the bias is a fixed, data-independent quantity.

\begin{proposition}[Linear growth of the prediction error]\label{P5.2}
Let $A_d\in\SpH(Q)$, let $\hat{z}_k = \Lift(x_k)$ be the lifted true trajectory and $r_k = \hat{z}_{k+1} -A_d\hat{z}_k$ the one-step residual. 
The model trajectory $z_{k+1}=A_dz_k$ with $z_0 = \hat{z}_0$ satisfies the inequality $\lVert z_k-\hat{z}_k \rVert_Q \le \sum_{j<k} \lVert r_j \rVert_Q$. 
If instead only $\lVert A_d \rVert_Q = g>1$ is known, the bound becomes $\sum_{j<k} g^{k-1-j} \lVert r_j \rVert_Q$.
\end{proposition}
\begin{proof}
With $e_k = z_k-\hat{z}_k$ one has $e_{k+1} = A_d e_k -r_k$ and $e_0 = 0$. Therefore, we have the inequality $\lVert e_{k+1} \rVert_Q \le \lVert A_d \rVert_Q \lVert e_k \rVert_Q + \lVert r_k \rVert_Q$.
\end{proof}

The bound concerns the lifted error and can be loose, but it is informative when the residuals are small. 
This explains why unconstrained models with a $Q$-gain above one can leave the data region: their bound grows geometrically.


\section{Choice of the Storage Metric}\label{S6}
The lifted pH triple $(\JK, \RK,Q)$ is determined by the data only up to congruence. 
Therefore, the storage metric $Q$ must be fixed by normalization.
The physically meaningful choice is \emph{Hamiltonian consistency}, which is defined as $\tfrac{1}{2} \Lift^\top Q \Lift = \Hamil$, making $\HK$ the physical energy. 
When the Hamiltonian $\Hamil$ is known, the admissible $Q$ form the solution set to an SDP feasibility problem.

For the damped pendulum with $\Hamil = \tfrac{1}{2} p^2 + \kappa -\cos\theta$, $\kappa>1$, and the \emph{dictionary} 
\begin{equation*}
    (1, \cos\theta, \sin\theta, p, p\cos\theta, p\sin\theta),
\end{equation*}
an explicit choice is given by 
\begin{equation}\label{E8}
Q = \begin{bmatrix}
2\kappa-1&-1&0\\ -1&1&0\\ 0&0&1
\end{bmatrix}
\oplus\diag(1-q,\,q,\,q), \quad\text{where}\quad 0<q<1.
\end{equation}
This matrix $Q$ is positive definite and satisfies $\tfrac{1}{2} \Lift^\top Q\Lift = \Hamil$ because $\cos^2\theta+\sin^2\theta=1$. 
The constant $\kappa$ is harmless, since only $\grad\Hamil$ enters \eqref{E1}.
If $\Hamil$ is unknown, a pragmatic alternative is to take $Q$ from a Lyapunov certificate of the unconstrained EDMD matrix \cite{math14214047}, 
or to alternate between the convex problem in $A_d$ for fixed $Q$ and an SDP in $Q$ for fixed $A_d$. 
Each subproblem is convex, but the alternation carries no global guarantee.

Two choices were tested. 
(C1) \emph{Energy-anchored}: coordinates $z=(z_E,z_X)$ with $H = \tfrac{1}{2} |z_E|^2$ and $Q = \mathrm{blkdiag}(I,\alpha I)$, so that $\tfrac{1}{2} z^\top Qz = H + \tfrac{\alpha}{2} |z_X|^2$. 
For the pendulum with $s=\sin(\theta/2)$, $H = \tfrac{1}{2} p^2 + 2 s^2$ is exactly quadratic in $(s,p)$, and $z_E = (p,2s)$. 
(C2) \emph{Data-adapted}: $Q_L$ is the solution of $A_c^\top Q_LA_c-Q_L=-I$ with $A_c = \widehat{A} / \max\bigl\{1,1.001 \rho(\widehat{A})\bigr\}$, rescaled such that $\tr Q_L=N$. 

Choice (C1) preserves the physical meaning of the storage, but it carries a large structural bias, because the lifted operator is non-normal and generally not a $Q$-contraction for such a $Q$. 
Choice (C2) makes $\widehat{A}$ almost a $Q_L$-contraction, so the bias is small.
However, in this case, the physical energy is not represented by the quantity $\tfrac{1}{2}z^\top Q_Lz$ and the guarantee concerns this fitted storage. 

Section~\ref{S7} quantifies this trade-off. 
The explicit $Q$ in \eqref{E8} was not benchmarked; the pendulum experiments below use the monomial dictionary in $(s,p)$.


\section{Numerical Examples}\label{S7}
We compare 
(a) EDMD \eqref{E2}, 
(a$'$) EDMD with the spectrum clipped to the unit disc, 
(b) restricted EDMD on $\SpH(Q)$, 
(c) restricted EDMD with Casimir constraints, and, for the pendulum, 
(d) restricted EDMD with the data-adapted metric $Q_L$ and 
(e) a continuous-time pH fit from finite-difference derivatives followed by a Cayley step. 
We solve the SDPs with CVXPY \cite{Diamond2016}, using Clarabel for $N\le 40$ and SCS beyond. 
If the returned matrix has a norm of $\lVert A_d\rVert_Q$ slightly above one, we divide it by this gain, except when Casimir constraints are active.
 In that case the residual is reported. 
Unless otherwise stated, we report medians over training sets and 30 held-out initial conditions. 
The error is the relative root mean square 
error of the original state over 50 steps.

\begin{example}[Pendulum]
We consider the pendulum equation $\ddot{\theta} = -\sin\theta -\gamma\dot{\theta}$, where $\gamma=0.02$, $\Delta t=0.1$, with observations $w = (\sin\tfrac{\theta}{2}, \dot{\theta})$, a dictionary of all monomials of degrees 1 to 5 in $w$ ($N=20$), and initial conditions sampled uniformly from $\theta\in[-2,2]$ and $\dot{\theta}\in[-1.4,1.4]$.
\textit{Settings:}
S1 (20 trajectories of 10 steps), S2 (100 of 15 steps), S3 (as S2 with noise of standard deviation $0.05$ on $w$), 10 training sets each.
The energy gain is $\max_t H(t)/H(0)$ over a free rollout of 1500 steps, computed on the energy coordinates.

Table~\ref{Tab1} and Figure~\ref{Fig1} show the safety side. 
The EDMD gains energy: in S2 the median peak is 54 times the initial energy (the worst case is about $10^4$). 
Clipping the spectrum eliminates divergence, but not energy growth (the worst case is 12 to 14 times the initial energy). 
For the energy-anchored contractions, the storage $\tfrac{1}{2} z^\top Qz$ cannot increase, therefore, the median energy gain is $1.00$.
However, the physical energy can still briefly increase in individual runs (the worst case is 1.7 in S1), because the extra coordinates $z_X$ return energy to $z_E$. 
The accuracy cost is significant, though:  with the physical metric, the 50-step error ranges from 0.12 to 0.31 in noise-free settings, compared to 0.03 to 0.05 for the EDMD. 
CT-pH achieves an accuracy of $\alpha=1$ in S1 and S2, consistent with both methods searching the same set up to closure. 
The data-adapted $Q_L$ is close to EDMD in S2.
In S3, EDMD, clipped EDMD, and $Q_L$ all reach the same noise-limited error of 0.224. 
However, the guarantee concerns the fitted storage, which has a physical energy gain of 1.02 and a worst case value of approximately 8 in S2.

\begin{table}[H]
\centering\footnotesize
\begin{minipage}[t]{0.48\textwidth}
\setlength{\tabcolsep}{.2cm}
\captionof{table}{Pendulum: median 50-step relative error (Err) and median energy gain (Gain). 
``div.'' indicates that the rollout diverged.
In S1, EDMD diverged in 89\,\% of the rollouts.}
\label{Tab1}
\centering
\adjustbox{max width=\linewidth}{
\begin{tabular}{@{}lcccccc@{}}
\toprule
 & \multicolumn{2}{c}{S1 (20$\times$10)}
 & \multicolumn{2}{c}{S2 (100$\times$15)}
 & \multicolumn{2}{c}{S3 (noisy)}\\
\cmidrule(lr){2-3}\cmidrule(lr){4-5}\cmidrule(lr){6-7}
Method & Err & Gain & Err & Gain & Err & Gain\\
\midrule
EDMD                    & 0.034 & div. & 0.048 & 54   & 0.224 & 1.01\\
EDMD, clipped           & 0.039 & 1.12 & 0.047 & 1.19 & 0.224 & 1.01\\
CT-pH + Cayley          & 0.246 & 1.00 & 0.314 & 1.00 & 0.332 & 1.00\\
$Q$-contr., $\alpha=1$   & 0.243 & 1.00 & 0.314 & 1.00 & 0.586 & 1.00\\
$Q$-contr., $\alpha=10$  & 0.117 & 1.00 & 0.179 & 1.00 & 0.540 & 1.00\\
$Q$-contr., $Q_L$        & 0.090 & 1.00 & 0.053 & 1.02 & 0.224 & 1.01\\
\bottomrule
\end{tabular}}
\end{minipage}
\hfill
\begin{minipage}[t]{0.48\textwidth}
\vspace{0pt}\centering
\includegraphics[width=\linewidth]{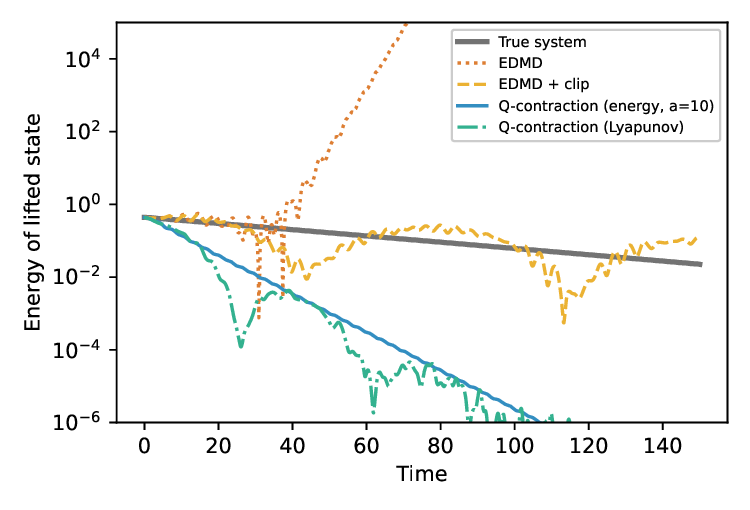}
\captionof{figure}{Energy of the lifted energy coordinates along one rollout (S1).
The energy-metric contraction exhibits monotonic decay (faster than the true dynamics, indicating structural bias), 
whereas the $Q_L$ model follows the true behavior more closely, although its physical energy is not monotone.
}\label{Fig1}
\end{minipage}
\end{table}



Figure~\ref{Fig2} varies the amplitude of the initial conditions (5 training sets per point). 
With the additional coordinates at their physical scale, the energy-anchored model achieves 0.6\,\% error at the smallest amplitude and 56\,\% at full amplitude.
Meanwhile, EDMD and the $Q_L$ model reach approximately 4\,\% and 4.6\,\%, respectively. 
This is consistent with the structural bias $\dist_{\widehat{G}}(\xi^\star,\mathcal{C})$ of Proposition~\ref{P5.1} growing with the nonlinearity of the data.

\begin{figure}[H]
\centering
\includegraphics[width=0.62\textwidth]{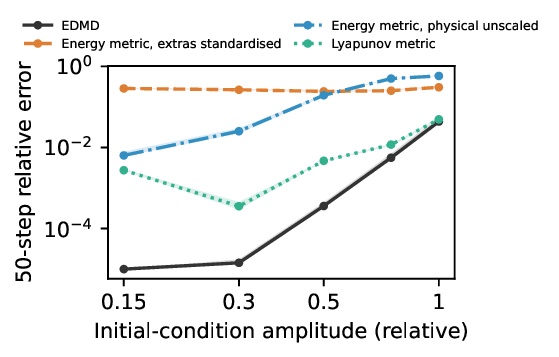}
\caption{Pendulum, 50-step relative error against the amplitude of the initial conditions (median, interquartile band; errors below $10^{-5}$ are drawn at $10^{-5}$).}
\label{Fig2}
\end{figure}
\end{example}

\begin{example}[Damped Rigid Body with a Casimir]
For the angular momentum $x\in\R^3$, $H(x) = \tfrac{1}{2} \sum_i x_i^2/I_i$ with $I=(1,2,3)$, and $J(x) = \hat{x}$, the function $C(x) = \tfrac{1}{2} \lVert x \rVert^2$ is a Casimir of $J$. 
Along a damped flow, $C$ is conserved only if also $R(x)\nabla C=0$ (the condition $c \in \ker J \cap \ker R$ in Theorem~\ref{T3.2}$(ii)$). 
With ordinary damping, the norm $\lVert x\rVert^2$ decays and the constraint $c^\top A_d = c^\top$ would be wrong. 
We therefore use the dissipation $R(x)=\gamma\,\hat{x}^\top\hat{x}$, i.e., $\dot{x} = x \times \nabla H +\gamma\,x \times(x \times \nabla H)$ with $\gamma=0.5$, for which $\dot{H} = -\gamma \lVert x \times \nabla H \rVert^2 \le0$ and $\dot{C}=0$. 
No torque is applied, since with an input $\dot{C}=x^\top Bu\ne0$, and $c^\top B_d = 0$ cannot hold for constant $B$. 
The dictionary consists of the monomials of degrees 1 to 3 in $x$ ($N=19$), which contain $H$ and $C$ exactly, and we use $Q = \mathrm{blkdiag}(I,10\,I)$ with the energy-anchored coordinates $z_E=(x_i/\sqrt{I_i})_i$. 

\textit{Settings:} R1 (30 trajectories of 10 steps), R2 (100 of 15), R3 (as R2, with noise of standard deviation $0.02$ on $x$); 5 training sets; rollouts of 500 steps.

\begin{table}[H]
\centering
\setlength{\tabcolsep}{0.4cm}
\caption{Rigid body: median 50-step error (Err), median energy gain $\max_t H(t)/H(0)$ (Gain), and median Casimir drift $\max_t |C(t)-C(0)|/C(0)$ (Cas.) over 500 steps.}\label{tab:rigid}
\centering\small
\adjustbox{max width=\textwidth}{
\begin{tabular}{@{}lccccccccc@{}}
\toprule
 & \multicolumn{3}{c}{R1 (30$\times$10)} & \multicolumn{3}{c}{R2 (100$\times$15)} & \multicolumn{3}{c}{R3 (noisy)}\\
\cmidrule(lr){2-4}\cmidrule(lr){5-7}\cmidrule(lr){8-10}
Method & Err & Gain & Cas. & Err & Gain & Cas. & Err & Gain & Cas.\\
\midrule
EDMD & 0.31 & $4{\cdot}10^{13}$ & $2{\cdot}10^{-5}$ & 0.22 & $8{\cdot}10^{7}$ & $8{\cdot}10^{-10}$ & 0.22 & $2{\cdot}10^{5}$ & 76\\
EDMD, clipped & 0.37 & 12.7 & $10^{-11}$ & 0.28 & 2.5 & $2{\cdot}10^{-12}$ & 0.26 & 1.17 & 0.71\\
$Q$-contr. & 0.35 & 1.00 & 1.08 & 0.32 & 1.00 & 1.10 & 0.37 & 1.00 & 1.05\\
$Q$-contr.\ + Casimir & 0.34 & 1.00 & $7{\cdot}10^{-9}$ & 0.33 & 1.00 & $2{\cdot}10^{-8}$ & 0.37 & 1.00 & $3{\cdot}10^{-8}$\\
\bottomrule
\end{tabular}}
\end{table}

Three observations follow from Table~\ref{tab:rigid}. 
First, the Casimir constraint matters only for noisy data. 
With exact data, the identity $c^\top Z_Y = c^\top Z_X$ holds, so the least-squares matrix already satisfies $c^\top \widehat{A} = c^\top$ when $Z_X$ has full row rank, and EDMD preserves $C$ (R1, R2).
What EDMD does not control is the energy (gains up to $10^{13}$; 41\,\% of the R1 rollouts diverged). 
With noise (R3), 
EDMD loses the Casimir completely (median drift 76), and clipping does not restore it. 

Second, a $Q$-contraction without the constraint does not preserve $C$ (drift above 100\,\%), since nothing in the contraction condition ties $C$ to its initial value. 

Third, the Casimir constraint costs no accuracy (0.34 versus 0.35 in R1) and reduces the drift to the order of $10^{-8}$. 
This is solver accuracy (Clarabel), not machine precision, and the contraction property holds with gain one up to the same tolerance. 
As for the pendulum, all constrained models remain less accurate than EDMD (errors 0.32 to 0.37 versus 0.22 in R2 and R3).
\end{example}

\begin{example}[Forced Chain with Ports]
Two unit masses between walls are coupled by hardening springs with forces $f_j = \delta_j + \delta_j^3$ in the elongations $\delta_j$ and damping $0.08$; a force input acts on the first mass and is held constant over each step of length $0.1$. 
The dictionary consists of the monomials of degrees 1 to 3 in $(x,p)$ ($N=34$). 
The supply is the work $u_k\,F(z_{k+1} - z_k)$, where the row vector $F$ extracts the displacement of mass~1 from $z$, so that $C_d = F(A_d - I)$ and $D_d = F B_d$, and \eqref{E7} is linear in $(A_d, B_d)$. 
Training uses 150 trajectories of 12 steps with random inputs in $[-0.5,0.5]$, without noise (F1) or with noise $0.02$ (F2), 5 training sets. 
EDMD violates \eqref{E7} strongly, 
whereas the smallest eigenvalue of the LMI matrix of the constrained fit is $-7.5\cdot10^{-5}$ (F1) and $-7.8\cdot10^{-5}$ (F2). 
This residual comes from the first-order solver, so passivity holds only numerically here. 
The 100-step error is 0.58 (F1) and 0.52 (F2) for the EDMD and 0.79 for the constrained model. 
Neither model is accurate, which indicates that this dictionary and data size are too small.
The free-decay rollouts of the unforced models did not diverge for any method. 
We report this example as a check that the port LMI can be imposed, not as evidence of a benefit.
\end{example}

\section{Conclusions}\label{S8}
For a fixed storage metric, the discrete-time Koopman matrices that admit a lifted pH realization are precisely the $Q$-contractions, which are an LMI-representable set.
Therefore, structure-preserving EDMD is a convex program. 
Its solution is passive and preserves the imposed linear Casimirs. 
The constraint cannot increase the estimation error when the population target is pH-compatible.
Experiments confirm these guarantees and demonstrate their cost.
 With physical energy as the storage metric, structural bias can dominate the error. 
However, a data-adapted metric removes most of the bias but loses the physical meaning of the storage.
Open questions include the consistency of the discrete-time estimator as $\Delta t\to 0$, the characterization of discrete-time bilinear pH maps, and the modeling error incurred by identifying a $Q$-contraction instead of the sampled Koopman matrix, which need not be a $Q$-contraction for finite $\Delta t$. 

\section*{Acknowledgments}
\noindent
M.~Ehrhardt acknowledges funding by the Deutsche Forschungsgemeinschaft (DFG, German Research Foundation), Project-ID 531152215, CRC~1701.

\section*{Declarations}

\subsection*{Data Availability} 
\noindent
Code and data are available from the authors upon request.

\subsection*{Conflict of Interest} 
\noindent
The authors declare that they have no known competing financial interests or personal relationships 
that could have appeared to influence the work reported in this letter.

\bibliographystyle{unsrt}
\bibliography{References}

\end{document}